\documentclass[a4paper, 11pt, reqno]{amsart}

\makeatletter
\usepackage{fullpage}

\usepackage{textcmds}  
\usepackage{amsmath, amssymb, amsfonts, amstext, verbatim, amsthm, mathrsfs, mathtools}
\usepackage[mathcal]{eucal}
\usepackage{stmaryrd}
\usepackage{microtype}
\usepackage{graphicx, psfrag, setspace, subfigure}
\usepackage[all,cmtip,2cell]{xy}
\UseTwocells
\usepackage{pgf,tikz,pgfplots}
\pgfplotsset{compat=1.15}
\usepackage{tikz-cd}
\usepackage{enumerate}
\usepackage{enumitem}
\usepackage{xspace}
\usepackage{cite}
\usepackage{pifont}
\usepackage{stmaryrd}
\usepackage{bbm}
\usepackage{fix-cm}
\usepackage{comment}
\usetikzlibrary{matrix,arrows,decorations.pathmorphing,decorations.pathreplacing}
\tikzset{commutative diagrams/diagrams={baseline=-2.5pt},commutative diagrams/arrow style=tikz}
\usepackage{float}
\usepackage{setspace} 

\usepackage[colorlinks=true,linkcolor=blue,citecolor=blue,urlcolor=blue,citebordercolor={0 0 1},urlbordercolor={0 0 1},linkbordercolor={0 0 1}]{hyperref} 
\usepackage[shortalphabetic]{amsrefs} 
\usepackage[nameinlink]{cleveref}

\def\makeCal#1{%
\expandafter\newcommand\csname c#1\endcsname{\mathcal{#1}}}
\def\makeBB#1{%
\expandafter\newcommand\csname b#1\endcsname{\mathbb{#1}}}
\def\makeFrak#1{%
\expandafter\newcommand\csname f#1\endcsname{\mathfrak{#1}}}

\count@=0
\loop
\advance\count@ 1
\edef\y{\@Alph\count@}%
\expandafter\makeCal\y
\expandafter\makeBB\y
\expandafter\makeFrak\y
\ifnum\count@<26
\repeat

\newcommand\Z{\mathbb Z}
\newcommand\C{\mathbb C}

\newcommand{\set}[1]{\left\{{#1}\right\}}

\newcommand\isoto{\stackrel{\sim}{\To}}
\newcommand\id{\mathrm 1}

\newcommand\To{\longrightarrow}

\newcommand\Hom{\operatorname{Hom}}
\newcommand\End{\operatorname{End}}

\newcommand\Ext{\operatorname{Ext}}

\newcommand\Ker{\operatorname{Ker}}

\renewcommand\P{\mathbb P}
\newcommand\Gr{\operatorname{Gr}}

\newcommand\Pic{\operatorname{Pic}}

\newcommand\rk{\operatorname{rank}}
\newcommand{\Tot}{\operatorname{Tot}}
\newcommand{\Sym}{\operatorname{Sym}}

\newcommand{\mat}[1]{\begin{pmatrix}#1\end{pmatrix}}
\newcommand{\smat}[1]{\left(\begin{smallmatrix}#1\end{smallmatrix}\right)}
\newcommand\MF{\operatorname{MF}}

\DeclareMathOperator{\DCoh}{D^b}

\DeclareMathOperator{\Tor}{Tor}

\newcommand*{\sheafhom}{\mathcal{H} \kern -.5pt om}

\newcommand{\beq}[1]{\begin{equation}\label{#1} }
\newcommand{\eeq}{\end{equation}}
\newcommand{\pgap}{\vspace{5pt}}

\theoremstyle{plain}
\newtheorem{prop}[equation]{Proposition}
\newtheorem{thm}[equation]{Theorem}
\newtheorem{lem}[equation]{Lemma}
\newtheorem{cor}[equation]{Corollary}

\theoremstyle{remark}
\newtheorem{rem}[equation]{Remark}

\theoremstyle{definition}

\newtheorem{notn}[equation]{Notation}

\makeatletter \@addtoreset{equation}{section} \makeatother

\newtheorem{keythm}{Theorem}

\let\oldtocsection=\tocsection
\let\oldtocsubsection=\tocsubsection
\let\oldtocsubsubsection=\tocsubsubsection
\renewcommand{\tocsection}[3]{\hspace{0em}\oldtocsection{#1}{#2}{#3}}
\renewcommand{\tocsubsection}[3]{ \hspace{1em} \oldtocsubsection{#1}{\small{#2}}{\small{#3}} }
\renewcommand{\tocsubsubsection}[3]{\hspace{2em}\oldtocsubsubsection{#1}{\small{#2}}{\small{#3}}}

\newcommand{\marginparstretch}{0.6}
\let\oldmarginpar\marginpar
\renewcommand\marginpar[1]{\-\oldmarginpar[\framebox{\setstretch{\marginparstretch}\begin{minipage}{\marginparwidth}{\raggedleft\scriptsize #1}\end{minipage}}]{\framebox{\setstretch{\marginparstretch}\begin{minipage}{\marginparwidth}{\raggedright\scriptsize #1}\end{minipage}}}}

\AtBeginDocument{%
   \def\MR#1{}
} 

\newcommand{\aand}{\quad\quad\mbox{and}\quad\quad}
\newcommand\ie{\emph{i.e.}~}

\makeatletter

\usepackage{babel}

\begin{document}

\title{The Fano of lines, the Kuznetsov component, and a flop}

\author{Kimoi Kemboi and Ed Segal}

\begin{abstract}The Kuznetsov component of the derived category of a cubic fourfold is a `non-commutative K3 surface'. Its symmetric square is hence a `non-commutative hyperk\"ahler fourfold'. We prove that this category is equivalent to the derived category of an actual hyperk\"ahler fourfold: the Fano of lines in the cubic. This verifies a conjecture of Galkin.

One of the key steps in our proof is a new derived equivalence for a specific 12-dimensional flop. 

\end{abstract}

\maketitle


\section{Introduction}

Let $Y\subset \P^5$ be a smooth cubic 4-fold. There is an associated Fano variety $F_Y$ of lines in $Y$, which is a smooth hyperk\"ahler 4-fold. More categorically, there is an associated Kuznetsov component 
$$\cA_Y\subset D^b(Y)$$
which famously behaves as a `non-commutative K3 surface'. For special classes of cubics this category $\cA_Y$ is actually equivalent to $D^b(S)$ for some genuine K3 surface $S$, and in these cases there is another hyperk\"ahler 4-fold available: the Hilbert scheme $S^{[2]}$.  So it is reasonable to suppose that these two hyperk\"ahler 4-folds might be related. 

Moreover, if we are only interested in categories, then there is no need to assume that $\cA_Y$ is geometric. If $S$ does exist then by \cite{BKR01} we have
$$D^b(S^{[2]}) \cong D^b_{\Z_2}(S\times S) \cong \Sym^2 \cA_Y$$
where the latter operation -- taking the symmetric square of a category -- is a purely formal operation. So for any $Y$ we can imagine there would be a relationship between $\Sym^2 \cA_Y$ and the derived category of $F_Y$. And indeed there is.

\begin{keythm}\label{mainthm1}
For any smooth cubic 4-fold $Y$ we have an equivalence:
 $$D^b(F_Y) \cong \Sym^2 \cA_Y $$
\end{keythm}

This result was apparently conjectured first by Galkin,\footnote{After completing this paper we were made aware of a 2017 talk by Galkin in which he not only states the conjecture but also explicitly anticipates the outline of our proof, see Remark \ref{galkin}.} and a weaker version of the statement was proven by Belmans-Fu-Raeschelders \cite{BFR22}. We learned about it from a recent paper of Bottini-Huybrechts \cite{HB25} who prove it for certain classes of cubics, see that paper for a nice outline of the history and also a discussion of the Hodge theory (which we will not touch on here). Some related conjectures are surveyed in \cite{BFM_2024}.\footnote{In the language of that paper, our result proves (\textbf{FM})$\implies$(\textbf{DF}).}
\pgap

We prove Theorem \ref{mainthm1} by using some well-established tricks in the world of matrix factorizations to reduce it to a simpler piece of higher-dimensional geometry. Let
$$X_0 \subset \Sym^3 (\C^6)^\vee$$
denote the space of cubic forms with a 4-dimensional kernel, \ie which can be factored through $\Sym^3 T^\vee$ for some 2-dimensional quotient $T$ of $\C^6$. This is a singular 12-fold. It has an obvious crepant resolution, which we call $X_-$, given by the total space of the vector bundle $\Sym^3 T^\vee$ over $\Gr(6,2)$. As we shall see, $X_0$ also has a second crepant resolution $X_+$, which is a smooth orbifold.

\begin{keythm}\label{mainthm2}
$X_+$ and $X_-$ are derived equivalent.
\end{keythm}
We explain the argument connecting Theorems \ref{mainthm1} and \ref{mainthm2} in Section \ref{BtoA}. It follows a similar pattern to various other examples, see for example \cite{ADS}*{Sect.~1.3}.

Theorem \ref{mainthm2} (Theorem \ref{tilting} below) is a special case of the well-known Bondal-Orlov-Kawamata conjecture that flops should induce derived equivalences. It may be of independent interest as it does not fall immediately to existing techniques, we have to do a bit of extra work. Thus it provides one more data point for the general conjecture. 
\pgap

\textbf{Acknowledgements.}
 We thank Daniel Huybrechts for encouraging correspondence and for pointing us to \cite{Gal}.

\section{From B to A}\label{BtoA}

Fix a non-singular cubic polynomial $f$ in 6 variables, so $Y=(f)\subset \P^5$ is a smooth cubic 4-fold. By Orlov's theorem \cites{Orlov_2009,S11} the Kuznetsov component $\cA_Y\subset D^b(Y)$ is equivalent to a category of matrix factorizations on an orbifold:\footnote{This statement requires an additional R-charge/grading on the orbifold which acts diagonally with weight $1/3$. This R-charge is present throughout our constructions but it plays no real role in the arguments so we leave it to the reader to insert it where needed.}
$$\cA_Y \;\cong\; \MF\!\big(\, [\C^6/\Z_3], \, f\,\big) $$
It follows that the symmetric square of $\cA_Y$ is given by matrix factorizations on the symmetric square of this orbifold, \ie
\beq{sym2}\Sym^2 \cA_Y \;\cong\; \MF\!\big(\, [(\C^6\times \C^6) /\Gamma ], \, f_1+f_2\,\big) \eeq
where $\Gamma$ is the semidirect product:
\beq{gamma}\Gamma =  (\Z_3\times \Z_3)\rtimes \Z_2 \eeq

Now consider the Fano variety of lines $F_Y$. This is a subvariety of $\Gr(2,6)$ cut out by a transverse section of $\Sym^3 U^\vee$, where $U$ is the tautological rank 2 bundle. By Kn\"orrer periodicity (\cite{Orlov_2004} \emph{etc.}), we have an equivalence
\beq{KP}D^b(F_Y) \;\cong\; \MF(X_-, W)\eeq
where $X_-$ is the total space of the dual bundle
$$\Sym^3 U \to \Gr(2,6) $$
and $W$ is the superpotential induced by $f$. We can build this space $X_-$ as a GIT quotient of the vector space
$$ \Hom(U, \C^6) \oplus \Sym^3 U$$
by the group $GL(U)\cong GL_2$. If we write $(\Phi,\psi)$ for the coordinates on the two summands, then the superpotential is:
$$W= (f\circ \Sym^3 \Phi)(\psi) $$
Note that this GIT problem is `Calabi-Yau' because $GL_2$ is acting trivially on the determinant of the vector space. Consequently $X_-$ is a (non-compact) Calabi-Yau, \ie it is a crepant resolution of the underlying quotient singularity. 
\pgap

Since it is more comfortable to think in terms of cubic polynomials rather than symmetric trivectors, we replace $U$ with $T=U^\vee$, so our GIT problem is the stack:
\beq{cX}\cX =  [\, \Hom(\C^6, T)\oplus \Sym^3 T^\vee \; /\; GL(T) \,] \eeq
This GIT problem has two possible stability conditions. For one choice, the semistable locus is $\{\rk \Phi =2\}$ and the GIT quotient is  $X_-$. We claim (see Section \ref{stability}) that for the other stability condition, semistability is the following two conditions:
\begin{itemize}\setlength{\itemsep}{3pt}
\item The cubic $\psi$ is non-zero with at least two distinct roots.
\item The image of $\Phi$ is not contained in a double root of $\psi$. 
\end{itemize}
This other GIT quotient, which we denote $X_+$, is a smooth Calabi-Yau orbifold. Thus $X_+$ and $X_-$ are two crepant resolutions of the quotient singularity $X_0$ underlying the stack $\cX$, and we have a 12-dimensional flop:
$$X_- \dashleftarrow\dashrightarrow X_+ $$
As explained in the introduction, the Bondal-Orlov-Kawamata conjecture predicts that this flop should induce a derived equivalence, and our Theorem \ref{mainthm2} verifies this. More precisely we will prove:

\begin{thm}\label{tilting}
There are tilting generators $E_-\in D^b(X_-)$ and $E_+\in D^b(X_+)$ with the same endomorphism algebra
$$A = \End_{X_-}\!(E_-) = \End_{X_+}\!(E_+)$$
and hence $D^b(X_-)\cong D^b(A) \cong D^b(X_+)$. 
\end{thm}

\begin{rem}

\begin{enumerate}
\item The object $E_-$ is the `obvious' tilting vector bundle coming from Kapranov's exceptional collection on the Grassmannian. But the object $E_+$ is not quite a vector bundle, and this is where some extra work is necessary.
\item The algebra $A$ is a non-commutative crepant resolution of $X_0$ and an example of the general theory of \cite{SVdB17}. 
\end{enumerate}\end{rem}

The superpotential $W$ on $\cX$ restricts to give a superpotential on both $X_-$ and $X_+$. We can also view it as a central element of the algebra $A$. Then the following corollary is standard (\emph{e.g.} \cite{S11}*{Lem.~3.6}. 

\begin{cor}\label{tiltingMF}
\[
\MF(X_-, W) \;\cong\; \MF(A, W) \;\cong\; \MF(X_+, W)
\]
\end{cor}

Now let $X_+^o\subset X_+$ be the open substack where the cubic $\psi$ has distinct roots. The group $GL_2$ acts transitively on such cubics, thus in this locus we may fix $\psi$ to be 
$$\psi = t_1^3 + t_2^3 $$
and we are left with a residual action of the stabilizer of this cubic polynomial, which is the finite group $\Gamma$ \eqref{gamma}. So
\beq{X+o}X_+^o = [\, \C^6 \times \C^6\;/\; \Gamma\, ] \eeq
and moreover the superpotential $W$ in this locus is the function $f_1+f_2$. This is precisely the LG model \eqref{sym2} that describes $\Sym^2 \cA_Y$. 
\pgap

 Finally we claim (see Section \ref{critical}) that the critical locus of $W$ on $X_+$ is entirely contained in the open set $X_+^o$, which implies that the restriction functor
\beq{restriction}\MF(X_+, W) \to \MF(X_+^o, W)\eeq
is an equivalence \cite{Orlov_2004}. So Corollary \ref{tiltingMF} provides the central link in a chain of equivalences
$$D^b(F_Y) \,\cong\, \MF(X_-, W) \,\cong\, \MF(X_+, W)\,\cong\,\MF(X_+^o, W) \,\cong\, \Sym^2 \cA_Y $$
and proves Galkin's conjecture. 

\begin{rem} The smoothness of $Y$ is required only for the step where we restrict to $X_-^o$, \emph{i.e.} to have that \eqref{restriction} is an equivalence. The other steps work for all cubics $f$. 
\end{rem}

\begin{rem}\label{galkin} In \cite{Gal}, Galkin points out the first and last of the above equivalences and observes that $X_-$ and $X_+^o$ are related by the GIT problem $\cX$. Our contribution is to understand the role played by $X_+$.
\end{rem}

\section{Supporting details} \label{S:geometric_constructions}

\subsection{Cubics} The space $\Sym^3 T^\vee$ of cubics in two variables has an obvious stratification 
$$\{0\} = \Sigma_0 \subset \Sigma_1 \subset \Sigma_2 \subset \Sigma_3= \Sym^3 T^\vee$$
by the number of roots of the cubic. Within each $\Sigma_i\setminus \Sigma_{i-1}$ the $GL_2$ action is transitive and it will be convenient at times to use the following `standard forms' for a cubic $\psi$:
$$ t_1^3, \quad\quad t_1^2t_2, \quad \quad t_1^3 + t_2^3$$
Note that $\Sigma_2$ is a divisor, cut out by the discriminant polynomial of a cubic. The stratum $\Sigma_1$ is a surface.

\subsection{Representation theory}\label{repthy} Recall that $\Gamma$ is the finite subgroup \eqref{gamma} of $GL(T)\cong GL_2(\C)$ that stabilizes the cubic $t_1^3+t_2^3$. For later use, we collect some elementary facts on the representation theory of $\Gamma$:
\begin{enumerate}\setlength{\itemsep}{3pt}    
\item[(i)] $\Gamma$ has nine irreducible representations, of which three are 2-dimensional and six are 1-dimensional.
\item[(ii)] The restriction of $\det(T)$ to $\Gamma$ generates the group of characters $\Pic(B\Gamma)\cong \Z_6$.
\item[(iii)]\label{repthyfact} The restriction of $\Sym^2 T\otimes (\det T)^{-1}$ to $\Gamma$ contains a 1-dimensional summand which is the restriction of $\det(T)^3$.
\end{enumerate}

\subsection{Stability}\label{stability} 

Here we analyse GIT stability on the stack $\cX$ \eqref{cX}. There are two stability conditions given by either a positive or negative power of the character $\det$ of $GL(T)$. 

As usual we look for destabilizing 1-parameter subgroups. The most obvious is the centre
$$\lambda_0 = \{ t\id_T,\, t\in \C^*\}$$
which destabilizes either the locus $\Phi=0$ or the locus $\psi=0$. For the remainder we fix a basis of $T$, and seek 1-parameter subgroups which fix more than just the origin. Up to conjugation and scaling there are exactly two.

\begin{itemize}\setlength{\itemsep}{5pt}

\item $\lambda_1 = (0,1)$. For the negative stability condition, this destabilizes points where $\rk \Phi <2$. For the positive stability condition, it destablizes points where $\psi\in \Sigma_1$, \ie cubics with a triple root. 

\item $\lambda_2=(-1,2)$. For the negative stability condition this destabilizes the set:
$$\{ \rk \Phi <2 \mbox{ and the image of } \Phi \mbox{ is a root of } \psi\}, $$
but these points are already destabilized by $\lambda_1$. For the positive stability condition, it destabilizes:
$$\{ \rk \Phi <2 \mbox{ and the image of } \Phi \mbox{ is a double root of } \psi\}. $$

\end{itemize}

\noindent Thus, for the negative stability condition the unstable locus is exactly $\{\rk \Phi< 2\}$. There are no strictly semistable points, and the GIT quotient is the vector bundle
\[
X_- = \mbox{Tot}\big\{\; \Sym^3 T^\vee \to \Gr(6,2)\; \big\}
\]
as claimed above. For the positive stability condition, we see that:
\begin{itemize}\setlength{\itemsep}{5pt}
\item Points with $\psi\in \Sigma_1$ are unstable.
\item Points with $\psi\in \Sigma_3\setminus \Sigma_2$ are stable.
\item Points with $\psi\in \Sigma_2\setminus \Sigma_1$ are unstable if and only if the image of $\Phi$ lies in the double root. 
\end{itemize}
The GIT quotient $X_+$ contains an open set $X_+^o$ where $\psi\notin \Sigma_2$, this is equivalent to $[\C^{12}/\Gamma]$ as noted above \eqref{X+o}. Now consider the divisor:
$$D = X_+ \setminus X_+^o $$
Within $D$ we can set $\psi=t_1^2t_2$, so the residual symmetry group is exactly $\lambda_2$. Semistability requires that the image of $\Phi$ does not lie in the double root $\{t_1=0\}$, \ie the top row of the matrix is non-zero. It follows that 
\beq{D}D \cong \mbox{Tot}\big\{\; \cO(-2)^{\oplus 6} \to \P^5\; \big\} \eeq
We observe that there are no semistable points with infinite isotropy, so $X_+$ is a smooth orbifold. 

\begin{rem}\label{discriminant} $D$ is cut out by the discriminant of the cubic $\psi$, which is a degree 4 polynomial in the coefficients of $\psi$. Since $\psi\in \Sym^3 T^\vee$ its discriminant lies in $\Sym^4(\Sym^3 T)$, specifically in the 1-dimensional summand $(\det T)^6$. Hence $\cO(D) = (\det T)^6$. 
\end{rem}

\subsection{The flop}\label{flop} 

We now make a few more observations on the geometry of the flop between $X_+$ and $X_-$. They are not all logically necessary for our argument, but we found them reassuring.
\pgap

The quotient singularity underlying $\cX$ is, as stated in the introduction, the subvariety
$$X_0 \subset \Sym^3(\C^6)^\vee$$
of cubic forms which have ``rank 2''. Obviously $X_0$ is singular at the origin, but this singularity is not isolated. There is a non-compact locus of singularities 
$$X_0^{sg}\subset X_0$$
 consisting of the cubics which are ``rank 1'',  \ie cubics that are a cube of a linear form. 

The spaces $X_+$ and $X_-$ are crepant resolutions of the Gorenstein quotient singularity $X_0$. We claim that away from the origin they are related by a family of standard toric flops.
\pgap

To see this, fix a generic point in $X_0^{sg}$, which we can take to be the cubic $x_1^3$. The deformations of this cubic that stay in $X_0^{sg}$ to first order are given by adding the terms $x_1^3$ and $x_1^2x_i$ for $i=2,.., 6$. So, a transverse slice to $X_0^{sg}$ is given by the family
$$x_1^3 + x_1Q(x_2,...,x_6) + C(x_2,..., x_6)$$
for a general quadratic polynomial  $Q$ and cubic $C$. Now consider the preimage of this family in $\cX$. We must have that the first basis vector is not in the kernel of $\Phi$, so we can take the first column of $\Phi$ to be $\smat{1\\0}$, so $\psi$ is of the form:
$$\psi = t_1^3 + \alpha t_1t_2^2 + \beta t_2^3 $$
The residual symmetry is the 1-parameter subgroup $\lambda_1 = (0,1)\subset GL_2$, which acts on the second row of $\Phi$ with weights 1 and on $\alpha, \beta$ with weights $-2$ and $-3$. In $X_-$, semistability forces the second row of $\Phi$ to be non-zero, and in $X_+$, semistability forces $(\alpha, \beta)\neq (0,0)$. So in this locus, our flop is a trivial family, over the base $\mathbb{A}^5$, of the standard orbifold flop:

$$ \mbox{Tot}\big\{\; \cO(-2)\oplus \cO(-3) \to \P^4 \; \big\} \quad \dashleftarrow\dashrightarrow \quad
\mbox{Tot}\big\{\; \cO(-1)^{\oplus 5} \to \P^1_{2:3} \; \big\}  $$
\vspace{-1pt}

Now we examine what happens at the origin in $X_0$. The fibre over the origin in $X_-$ is the Grassmannian $\Gr(6,2)$. On the positive side, let us denote the fibre over the origin by:
$$F \subset X_+ $$
Then $F$ is the locus in $X_+$ where the image of $\Phi$ lies in a root of $\psi$. It's not immediately obvious what this space is, but we can understand it by intersecting with the open set $X_+^o$ and the complementary divisor $D$  \eqref{D}. The intersection of $F$ with $D$ is just the zero section:
$$F\cap D = \P^5 $$
The intersection of $F$ with $X_+^o$ is the subvariety in $[\Hom(\C^6, \C^2) /\Gamma]$ where the image of $\Phi$ lies in one of the three roots of the cubic $t_1^3+t_2^3$. This is three 6-dimensional linear subspaces meeting at the origin. So $F$ is a singular orbifold, \ie the quotient of a singular variety by a finite group. 

\begin{rem}\label{eulerchar}
Since $\Gamma$ has 9 irreducible representations we have an equality of Euler characteristics
$$\chi(X_+) \,=\, \chi(F) \,=\, \chi(B\Gamma) + \chi(\P^5) \,=\, 9 + 6 \,=\, 15 \,=\, \chi(\Gr(2,6))  \,=\, \chi(X_-)$$
which is a necessary condition for derived equivalence. 
\end{rem}

Finally, we observe that $F$ has a resolution by a weighted projective space:
\beq{tildeF}\widetilde{F} \cong \P^6_{1:...:1:6} \; \To \;F \eeq
We construct this resolution by making an additional choice of a line $L\subset T$ that is a root of $\psi$, and which contains the image of $\Phi$. Away from $\Phi=0$, this is no extra data, but at the singular point it separates the three branches whilst also reducing the isotropy group from $\Gamma$ to $\Z_6$. To see that the resulting resolution $\widetilde{F}$ really is the claimed weighted projective space, we set $L=\{t_2=0\}$ and choose $\psi$ to be of the form:
$$\psi = t_2( t_1^2 + \gamma t_2^2)$$
The image of $\Phi$ lies in $\{t_2=0\}$, \ie the bottom row of $\Phi$ is zero, and the residual symmetry is the 1-parameter subgroup $\lambda_2$. It acts on the top row of $\Phi$ with weight $-1$ and on the coefficient $\gamma$ with weight $-6$.

\subsection{The critical locus}\label{critical} 

Here we analyse the critical locus of the superpotential $W$ on the stack $\cX$. Recall that 
$$W =(\psi\circ \Sym^3 \Phi)(f) $$
where $\psi\in \Sym^3 T^\vee$, $\Phi\in \Hom(\C^6, T)$, and $f$ is a fixed cubic in 6 variables, the defining equation of our cubic fourfold $Y$.  
\pgap

First, observe that $W$ is linear in $\psi$, so $\partial_\psi W=0$ means $\Sym^3 \Phi(f)=0$. This says that the pullback of the cubic $f$ via $\Phi$ is the zero cubic in two variables. In particular, either

\begin{itemize}\setlength{\itemsep}{5pt}
\item $\Phi=0$; or 
\item $\rk(\Phi)=1$ and the image of $\Sym^3 \Phi$ is a point in $Y$; or 
\item $\rk(\Phi)=2$ and the image of $\Phi$ is a point of $F_Y$. 
\end{itemize}

Now we look at the other derivatives.  We use the stratification $\Sigma_i$ of $\Sym^3 T^\vee$ introduced above and deal with each stratum separately. 
\pgap

\begin{itemize}\setlength{\itemsep}{7pt}

\item[($\Sigma_0$)] If $\psi=0$, then $W_\psi\equiv 0$, so the $\Phi$ derivatives vanish.

\item[($\Sigma_1$)] If $\psi=t_1^3$, then $W_\psi = f(\Phi_{11}, ..., \Phi_{16})$. Since $f$ is non-singular, the $\Phi$ derivatives vanish exactly at $\Phi_{1\bullet}=0$, \ie when $\rk{\Phi}\leq 1$ and the image of $\Phi$ is in the root of $\psi$. 

\item[($\Sigma_3$)] If $\psi=t_1^3 + t_2^3$, then $W_\psi = f(\Phi_{1\bullet}) + f(\Phi_{2\bullet})$, and the $\Phi$ derivatives vanish only at $\Phi=0$. 

\item[($\Sigma_2$)] Set $\psi=t_1^2t_2$. Write $\widehat{f}$ for the symmetric trilinear map corresponding to the cubic form $f$, so
$$W_\psi = \widehat{f}(\Phi_{1\bullet}, \Phi_{1\bullet}, \Phi_{2\bullet}) $$
This is linear in the variables $\Phi_{2i}$, and we have:
$$\partial_{\Phi_{2i}} W_\psi = \partial_i f (\Phi_{11},..., \Phi_{16}) $$
Again these vanish exactly when $\Phi_{1\bullet}=0$, \ie when $\rk{\Phi}\leq 1$ and the image of $\Phi$ is in the double root of $\psi$.  And since $W_\psi$ is quadratic in the $\Phi_{1i}$ variables, the derivatives $\partial_{\Phi_{1i}} W_\psi$ also vanish at these points. 
\end{itemize}

Now we compare this with our description of the semistable loci from Section \ref{stability}. In $X_-$ we have $\rk \Phi=2$, so the only critical points are along
$$F_Y \subset \Gr(6,2) \subset X_- $$
This is a standard part of the Kn\"orrer periodicity equivalence \eqref{KP}. More importantly, we see that in $X_+$, there is a single critical point 
$$\{\Phi=0\} \in X_+^o$$
since all the critical points with $\psi\in \Sigma_2$ are unstable. This is essential for the restriction functor \eqref{restriction} to be an equivalence.

\section{The proof of Theorem B}

We will now prove Theorem \ref{tilting} by constructing tilting generators of $\DCoh(X_-)$ and $\DCoh(X_+)$ with the same endomorphism algebra. We prove the two sides separately as Propositions \ref{tiltingX-} and \ref{tiltingX+} below.

We will use the techniques of \cite{SVdB17} and the machinery of grade restriction subcategories \cite{HL15}, \emph{a.k.a.} windows, so most of our arguments take place on the linear stack $\cX$. This approach helps to motivate our construction of the generator $E_+$ but it is not essential; the proofs could be written purely on $X_+$ and $X_-$. 

\subsection{Grade restriction subcategories}\label{sec.GR}

Recall that our GIT problem is the vector space
$$V = \Hom(\C^6, T)\oplus \Sym^3 T^\vee$$
with the action of the group:
$$G = GL(T) \cong GL_2 $$
As before, we write $\cX=[V/G]$ for the Artin quotient stack. We have two GIT stability conditions and hence two GIT quotients, which are open substacks $X_{\pm}=[V^{ss}_{\pm}/G]\subset \cX$. The theory developed in \cite{HL15} provides us with two subcategories
$$\cG_{\pm} \subset D^b(\cX)$$
that are equivalent to $D^b(X_\pm)$ via the natural restriction functors. 
\pgap

Recall (Section \ref{stability}) that for the negative stability condition, unstable points are, up to the action of $G$, destabilized by one of the two 1-parameter subgroups $\lambda_0, \lambda_1$.
This provides a \emph{KN stratification} of the unstable locus. For each $i$ we have a fixed subspace
$$Z_i = V^{\lambda_i}$$
and also an attracting subspace $V^{\lambda_i \geq 0}$ of points destabilized by $\lambda_i$, \ie the sum of the eigenspaces for non-negative eigenvalues of $\lambda_i$. There is an obvious stratification of the unstable locus into locally closed pieces
$$GV^{\lambda_0\geq 0} = V^{\lambda_0\geq 0}=\set{\Phi = 0} \aand G V^{\lambda_1\geq 0}\setminus V^{\lambda_0\geq 0} = \set{\rk \Phi = 1}$$
but we also need to be precise about how this stratification interacts with the fixed and attracting subspaces.

\newcommand{\notate}[6]{\smat{\vec{#1}\\ \vec{#2}}\!, (#3, #4, #5, #6)}

\begin{notn} \label{N:eigenspace_decomposition}
Since we have fixed a basis of $T$, a point in $V$ consists of a matrix $\Phi$ and a 4-dimensional vector $\psi$, which we will write as
$$\notate{\Phi_1}{\Phi_2}{\psi_1}{\psi_2}{\psi_3}{\psi_4} $$
where the components $\vec{\Phi_i}$ are 6-vectors, the rows of $\Phi$. In such a presentation, we will use $\square$ to denote a component that is nonzero, and $\star$ to denote a component that is either zero or nonzero.
\end{notn}

With this notation in place, we introduce the locally closed subvarieties:

$$Z^-_0 = Z_0 =\set{\notate{0}{0}{0}{0}{0}{0}}\hspace{4cm} Y^-_0 = V^{\lambda_0 \geq 0} =\set{\notate{0}{0}{\star}{\star}{\star}{\star}}$$
$$Z^-_1 = Z_1\setminus Y^-_0 = \set{\notate{\square}{0}{\star}{0}{0}{0}}\hspace{3cm} Y^-_1 = V^{\lambda_1 \geq 0} \setminus Y^-_0= \set{\notate{\square}{0}{\star}{\star}{\star}{\star}}$$
The unstable locus on the negative side is $Y^-_0\sqcup GY^-_1$, and both pieces are smooth.  
\pgap

For the positive stability condition, we have a similar stratification of the unstable locus, except that (i) we have a third destabilizing 1-parameter subgroup $\lambda_2$ to consider, so the stratification has three pieces, and (ii) the destabilized points are now $V^{\lambda_i \leq 0}$ instead of $V^{\lambda_i \geq 0}$. Using the same notation as above, we have:

$$    Z_0^+ = Z_0= \set{ \notate{0}{0}{0}{0}{0}{0}}\hspace{3cm} Y_0^+ = V^{\lambda_0 \leq 0}= \set{ \notate{\star}{\star}{0}{0}{0}{0} } $$

$$    Z_1^+ = Z_1\setminus Y_0^+ = \set{ \notate{\star}{0}{\square}{0}{0}{0}} \hspace{3cm} Y_1^+ = V^{\lambda_1\leq 0}\setminus Y_0^+= \set{ \notate{\star}{\star}{\square}{0}{0}{0} }$$

$$    Z_2^+ = Z_2 \setminus GY_1^+ =  \set{\notate{0}{0}{0}{\square}{0}{0}} \hspace{2cm} Y_2^+ = V^{\lambda_2 \leq 0} \setminus GY_1^+ = \set{\notate{0}{\star}{\star}{\square}{0}{0} } $$
The unstable locus is $Y_0^+ \sqcup GY_1^+ \sqcup GY_2^+$, matching the description we gave in Section \ref{stability}.

We also need the following integer invariants:
$$ \eta_i = \sum \mbox{positive eigenvalues of }\lambda_i\mbox{ on }V
\;-\; \sum \mbox{positive eigenvalues of }\lambda_i\mbox{ on }\mathfrak{g} $$
Note that because $\cX$ satisfies the Calabi-Yau condition, replacing $\lambda_i$ with $-\lambda_i$ in this definition gives the same answer. In terms of the statification, $\eta_i$ computes the total $\lambda_i$-weight of the conormal bundle $N^\vee_{GY^{\pm}_i}V$. In our example, it is easy to compute:
$$\eta_0 = 12\quad\quad\quad \eta_1 = 6-1=5 \quad \quad\quad \eta_2 = 15 - 3 =12 $$

Given an object $\cF\in D^b(\cX)$, the homology sheaves $$\cH^\bullet(\cF\vert_{Z_i^{\pm}})$$ of the restriction of $\cF$ to $Z_i^{\pm}$ split into weight spaces for the $\lambda_i$ action. We define the \emph{grade restriction subcategories} $\cG_{\pm} \subset D^b(\cX)$ to be the full subcategories of objects $\cF$ that obey the following \emph{grade restriction rules}:
\beq{GRRs}\mbox{The }\lambda_i\mbox{-weights of }\cH^\bullet(\cF\vert_{Z_i^{\pm}})\mbox{ lie in }\; \big[ -\!\tfrac12\eta_i, \,\tfrac12\eta_i\big) \eeq
Here $i$ runs over the pieces of the stratification, so $\cG_-$ is defined by two grade restriction rules, and $\cG_+$ by three. 

\begin{thm}\label{HL}\cite{HL15} The restriction functors $D^b(\cX) \to D^b(X_\pm)$ induce equivalences:
$$\cG_{\pm} \isoto D^b(X_\pm)$$
\end{thm}

In fact, there are many possible grade-restriction subcategories, since for each $i$ one can shift the interval in \eqref{GRRs} by an integer and the above theorem still holds. For our purposes, it is simplest to center the intervals on $0$.
\pgap

In general, it can be hard to determine whether an object lies in a grade restriction category, but there is a class of objects for which it is easy. Pick an irreducible representation of $G$, indexed by some dominant weight $\chi= (a,b)$. There is a corresponding vector bundle on the stack $\cX$ which we denote by
$$\cT_\chi := \Sym^{a-b} T \otimes (\det T)^b $$
To see if $\cT_\chi$ satisfies \eqref{GRRs} we just pair $\lambda_i$ with both $\chi$ and its Weyl partner $\chi'=(b,a)$ and see if the answers lie in the right interval.

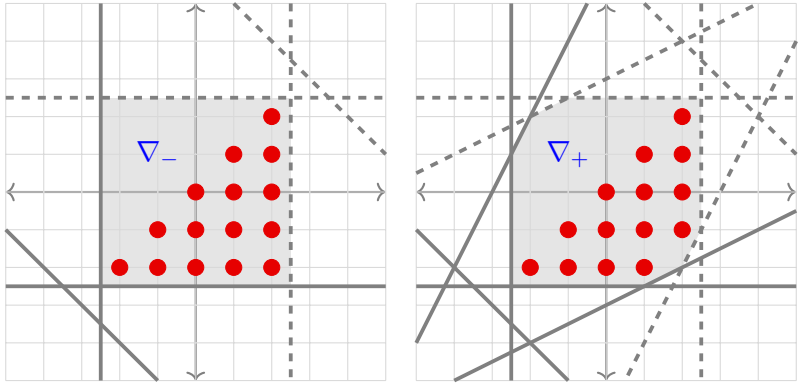
\begin{figure}[ht]
\begin{tabular}{ll}

\begin{tikzpicture}[scale=1]

\fill[fill=gray,fill opacity=0.2]   (-1.25,-1.25) -- (-1.25,1.25) -- (1.25,1.25) -- (1.25,-1.25) -- cycle;

\draw[<->, thick, color=gray] (-2.5,0) -- (2.5,0);
\draw[<->, thick, color=gray] (0,-2.5) -- (0,2.5);
\draw[step=0.5, thin, color=gray!30!white] (-2.5,-2.5) grid (2.5,2.5);

\draw[dashed, line width=0.5mm, color=gray] (0.5,2.5) -- (2.5,0.5);
\draw[dashed, line width=0.5mm, color=gray] (1.25,2.5) -- (1.25,-2.5);
\draw[dashed, line width=0.5mm, color=gray] (-2.5,1.25) -- (2.5,1.25);

\draw[-, line width=0.5mm, color=gray] (-2.5,-0.5) -- (-0.5,-2.5);
\draw[-, line width=0.5mm, color=gray] (-1.25,-2.5) -- (-1.25,2.5);
\draw[-, line width=0.5mm, color=gray] (-2.5,-1.25) -- (2.5,-1.25);

\node[blue] at (-0.5,0.5) {$\nabla_-$};

\foreach \i in {  (-1,-1), (-0.5,-0.5), (0,0), (0.5,.5),(1,1), (-0.5,-1),(0,-0.5),(0.5,0),(1,0.5), (0,-1), (0.5,-0.5), (1,0), (0.5,-1), (1,-0.5), (1,-1)}
{ \filldraw[red!90!black] \i circle (3 pt); }

\end{tikzpicture}

&

\begin{tikzpicture}[scale=1]

\fill[fill=gray,fill opacity=0.2]   (-1.25,-1.25) -- (-1.25,0.5) -- (-1,1) -- (-.5,1.25) -- (1.25,1.25) -- (1.25,-0.5) -- (1,-1) -- (0.5,-1.25) -- cycle;

\draw[<->, thick, color=gray] (-2.5,0) -- (2.5,0);
\draw[<->, thick, color=gray] (0,-2.5) -- (0,2.5);
\draw[step=0.5, thin, color=gray!30!white] (-2.5,-2.5) grid (2.5,2.5);

\draw[dashed, line width=0.5mm, color=gray] (0.5,2.5) -- (2.5,0.5);
\draw[dashed, line width=0.5mm, color=gray] (1.25,2.5) -- (1.25,-2.5);
\draw[dashed, line width=0.5mm, color=gray] (-2.5,1.25) -- (2.5,1.25);
\draw[dashed, line width=0.5mm, color=gray] (2.5,2) -- (0.25,-2.5);
\draw[dashed, line width=0.5mm, color=gray] (-2.5,0.25) -- (2,2.5);

\draw[-, line width=0.5mm, color=gray] (-2.5,-0.5) -- (-0.5,-2.5);
\draw[-, line width=0.5mm, color=gray] (-1.25,-2.5) -- (-1.25,2.5);
\draw[-, line width=0.5mm, color=gray] (-2.5,-1.25) -- (2.5,-1.25);
\draw[-, line width=0.5mm, color=gray] (-2.5,-2) -- (-.25,2.5);
\draw[-, line width=0.5mm, color=gray] (-2,-2.5) -- (2.5,-.25);

\node[blue] at (-0.5,0.5) {\textbf{$\nabla_+$}};

\foreach \i in {  (-1,-1), (-0.5,-0.5), (0,0), (0.5,.5),(1,1), (-0.5,-1),(0,-0.5),(0.5,0),(1,0.5), (0,-1), (0.5,-0.5), (1,0), (0.5,-1), (1,-0.5)}
{ \filldraw[red!90!black] \i circle (3 pt); }

\end{tikzpicture}

\end{tabular}
\caption{\footnotesize The red dots are the dominant weights corresponding to vector bundles that satisfy the grade restriction rules for $\cG_-$ (left) and $\cG_+$ (right).}\label{fig:windows_for_running_git}
\end{figure}

For $\cG_-$ this gives four inequalities cutting out a set $\nabla_-$ of 15 weights, and for $\cG_+$ it gives six inequalities cutting out a set $\nabla_+$ of 14 weights. See Figure \ref{fig:windows_for_running_git}. 

Thus, we have a set of 15 vector bundles:
\beq{15vbs}\left\{ \cT_\chi\;,\; \chi\in \nabla_-\right\} = \left\{ \Sym^c T\otimes (\det T)^b\; ,\;  c\leq 4,\ \; -2\leq b \leq 2-c\right\}\eeq
on $\cX$ which all lie in $\cG_-$. We shall prove in the next section that the direct sum of the restriction of these bundles to $X_-$ give a tilting generator of $D^b(X_-)$. Note that this is just the pullback to $X_-$ of Kapranov's exceptional collection on $\Gr(6,2)$ \cite{K88}.
\pgap

Unfortunately, things are not so simple on $X_+$. The bundle
$$\cT_{(2,-2)}=\Sym^4T\otimes(\det T)^{-2}$$ 
fails the grade restriction rule for $\lambda_2$, so it does not lie in $\cG_+$; and indeed one can show that the higher self-Exts of $\cT_{(2,-2)}\vert_{X_+}$ do not vanish (Remark \ref{selfExt}). Moreover, since the Euler characteristic of $X_+$ is 15 (Remark \ref{eulerchar}), the bundles from $\nabla_+$ alone cannot generate $D^b(X_+)$. So we will have to work a bit harder to find a tilting generator on the positive side.

\subsection{Weyman's complexes}

Consider the locus $GV^{\lambda_1 \geq 0}=\{\rk \Psi \leq 1 \}$, the closure of $GY_1^-$. This is a singular subvariety of $V$, and it has a Springer-type resolution:
$$r:G\times_P V^{\lambda_1 \geq 0} \To GV^{\lambda_1\geq 0} $$
Here, $P\subset G$ is the parabolic associated to $\lambda_1$, so $G/P=\P T\cong \P^1$. More explicitly, the resolution is the vector bundle:
$$G\times_P V^{\lambda_1 \geq 0}\; \cong \;\Tot\big\{\Hom(\C^6, L)\oplus \Sym^3 T^\vee \To \P T \big\} $$
where $L$ is the tautological line bundle on $\P T$. 

For any positive line bundle on this resolution, we have an associated sheaf $r_*(L^{-k})$ on $\cX$, supported on $GV^{\lambda_1\geq 0}$. Weyman \cite{W03} introduced a method for constructing free resolutions of such sheaves: we embed the resolution into $V\times \P T$, take the Koszul resolution of our sheaf, and push down. This produces a free resolution
$$C_{\lambda_1, k} \isoto r_*(L^{-k})$$
whose terms are sums of bundles $\cT_\chi$, and the weights that appear are easily calculated by Borel-Weil-Bott. This method works for any 1-parameter subgroup $\lambda\subset G$, and indeed more generally for other representations of other reductive groups.\footnote{For larger $G$ the resolution will be a vector bundle over some partial flag variety $G/P$.}

These complexes were used to great effect by \v{S}penko-Van den Bergh to construct non-commutative resolutions of quotient singularities \cite{SVdB17}. They are also relevant to GIT, since when we restrict $C_{\lambda_1, k}$ to $X_-$ we get an exact sequence of vector bundles. In the right circumstances such sequences can be used to show the derived category of the quotient is generated by vector bundles, \emph{e.g.} \cite{HLS20}. We can also view this as the process of mutating into the grade restriction category. 
\pgap

Now recall our set of 15 vector bundles \eqref{15vbs} that lie in $\cG_-$. We denote their sum by:
$$\cE_- = \bigoplus_{\chi\in \nabla_-} \cT_{\chi}$$

\begin{prop}\label{tiltingX-} The vector bundle $E_- =\cE_-\vert_{X_-} $
is tilting and its summands generate $D^b(X_-)$.
\end{prop}

\begin{proof}
From Theorem \ref{HL} we know:
$$\Ext^\bullet_{X_-}\!(E_-,E_-) = \Ext^\bullet_{\cX}(\cE_-, \cE_-) = \Ext^0_\cX(\cE_-, \cE_-) $$
On the Grassmannian itself, we know from \cite{K88} that every object can be resolved using the summands of $E_-$. It is a straight forward exercise to reprove this using the methods of \cite{SVdB17}; the complexes $C_{\lambda_1, k}$ produce enough exact sequences to reduce everything to $\nabla_-$. The exact same proof also works for $X_-$ since the $\Sym^3 T^\vee$ factor will have no effect on these constructions. 
\end{proof}

As before, we denote the endomorphism algebra of $E_-$ by
$$ A = \End_{X_-}\!(E_-)$$
Then $A$ is a non-commutative crepant resolution of $X_0$.

\subsection{The tilting sheaf on \texorpdfstring{$X_+$}{X+}}

Recall the subset of $\cX$ that is destabilized by $\lambda_2$ on the negative side:
$$GV^{\lambda_2\geq 0} = \{ \rk \Phi < 2 \mbox{ and the image of } \Phi \mbox{ is a root of } \psi\}$$
(This did not appear in our KN stratification since such points are already destabilized by $\lambda_1$). It is a singular subvariety of $V$, and it has a resolution
$$q:G\times_P V^{\lambda_2 \geq 0} \To GV^{\lambda_2\geq 0}  $$
Here, $P\subset G$ is the parabolic associated to $\lambda_2$, so $G/P=\P T\cong \P^1$, and explicitly:
\beq{resolution}G\times_P V^{\lambda_2 \geq 0}\; \cong \;\Tot\big\{\Hom(\C^6, L)\oplus \Sym^2 T^\vee\otimes (T/L)^\vee \To \P T \big\}  \eeq
If we intersect with $X_+$ then $GV^{\lambda_2 \geq 0}$ is precisely the subset $F$ discussed in Section \ref{flop}, the fibre of $X_+$ over the origin in $X_0$. Moreover, this resolution is the same as our previous $\widetilde{F}$ \eqref{tildeF}.
\pgap

Now consider the torsion sheaf $q_*(L^{-4})$ on $\cX$. By adjunction, there is a canonical map:
$$\Sym^4 T^\vee \to q_*(L^{-4}) $$
This is the final term of Weyman's free resolution. In particular, it is a surjection. We twist by a line bundle and take its kernel:
$$ \cK  = \Ker \Big(\; \cT_{(2,-2)}  \to q_*(L^{-4}) \otimes (\det T)^2 \Big)$$
This, it turns out, is the correct modification of $\cT_{(2,-2)}$. Notice that the restriction of $\cK$ to $X_-$ agrees with the bundle $\cT_{(2,-2)}$, but the restriction of $\cK$ to $X_+$ does not.

\begin{prop}\label{tiltingX+} Let
$$\cE_+ = \cK \oplus \bigoplus_{\chi\in \nabla_+} \cT_\chi$$
and let $E_+ = \cE_+\vert_{X_+}$. Then $E_+$ is a tilting sheaf, its summands generate $D^b(X_+)$, and it has endomorphism algebra:
$$\End_{X_+}\!(E_+) = A $$
\end{prop}
\begin{proof}
We will prove below the following claims: 
\begin{enumerate}\setlength{\itemsep}{5pt}
\item $\cK$ lies in the grade restricted subcategory $\cG_+$ (Lemma \ref{GRRforK}).
\item We have 
$$\End(\cE_+) = \End(\cE_-)$$
and $\Ext^k(\cE_+, \cE_+)=0$ for $k>0$ (Lemma \ref{EndK}). 
\item The sheaf $E_+$ spans $D^b(X_+)$ (Proposition \ref{E+spans}). 
\end{enumerate}
\pgap

From (1) it follows that $\Ext^\bullet_{X_+}\!(E_+, E_+) = \Ext^\bullet_\cX(\cE_+, \cE_+)$, and from (2) this is the algebra $A$. So the summands of $E_+$ generate a subcategory of $D^b(X_+)$ which is equivalent to $D^b(A)$ and hence admissible (since $A$ is smooth). By (3) this subcategory is the whole of $D^b(X_+).$
\end{proof}

It follows from above that the summands of $\cE_+$ generate $\cG_+$, and that the restriction functor from $\cG_+$ to $D^b(X_-)$ is an equivalence. In particular, our derived equivalence between $X_+$ and $X_-$ can also be expressed as the composition:
$$D^b(X_+)\; \stackrel{\sim}{\longleftarrow} \; \cG_+ \; \isoto\; D^b(X_-)$$
Thus $\cG_+$ provides the correct `window' for this derived equivalence. 

Note that we cannot replace $\cG_+$ with $\cG_-$ here because the restriction functor from $\cG_-$ to $D^b(X_+)$ is not an equivalence (Remark \ref{selfExt}). 

\subsection{More on \texorpdfstring{$\cK$}{\cK}} 
\begin{figure}

\begin{tikzpicture}[scale=1]

\fill[fill=gray,fill opacity=0.2]   (-1.25,-1.25) -- (-1.25,0.5) -- (-1,1) -- (-.5,1.25) -- (1.25,1.25) -- (1.25,-0.5) -- (1,-1) -- (0.5,-1.25) -- cycle;

\draw[<->, thick, color=gray] (-2.5,0) -- (2.5,0);
\draw[<->, thick, color=gray] (0,-2.5) -- (0,2.5);
\draw[step=0.5, thin, color=gray!30!white] (-2.5,-2.5) grid (2.5,2.5);

\draw[dashed, line width=0.5mm, color=gray] (0.5,2.5) -- (2.5,0.5);
\draw[dashed, line width=0.5mm, color=gray] (1.25,2.5) -- (1.25,-2.5);
\draw[dashed, line width=0.5mm, color=gray] (-2.5,1.25) -- (2.5,1.25);
\draw[dashed, line width=0.5mm, color=gray] (2.5,2) -- (0.25,-2.5);
\draw[dashed, line width=0.5mm, color=gray] (-2.5,0.25) -- (2,2.5);

\draw[-, line width=0.5mm, color=gray] (-2.5,-0.5) -- (-0.5,-2.5);
\draw[-, line width=0.5mm, color=gray] (-1.25,-2.5) -- (-1.25,2.5);
\draw[-, line width=0.5mm, color=gray] (-2.5,-1.25) -- (2.5,-1.25);
\draw[-, line width=0.5mm, color=gray] (-2.5,-2) -- (-.25,2.5);
\draw[-, line width=0.5mm, color=gray] (-2,-2.5) -- (2.5,-.25);

\foreach \i in {  (0,-1.5), (0,-1), (0,-.5), (0,0),(.5,.5), (1, .5),
(-1.5, -1.5), (-1,-1), (-.5, -1), (0,-1), (.5, -1)  }
{ \filldraw[green!60!black] \i circle (3 pt); }

\foreach \i in {  (1,-1)  }
{ \filldraw[blue] \i circle (4 pt); }

\foreach \i in {  (0, -1.5)  }
{ \draw[line width=0.5mm, red!90!black] \i circle (5 pt); }

\node[blue] at (2,-1) {$(2,-2)$};
\node[red] at (-.1,-2) {$(0,-3)$};
\end{tikzpicture}

\caption{\footnotesize The bundles appearing in the free resolution of $\cK$ (in green) together with $\cT_{(2,-2)}$ (in blue). The final term of the resolution is circled in red. The grade restriction rules for $\cG_+$ are also shown.} \label{fig:K}

\end{figure}
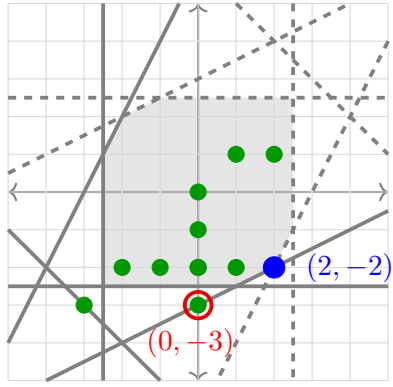

Using Weyman's method, the torsion sheaf $q_*(L^{-4})\otimes (\det T)^2$ has a free resolution by bundles of the form $\cT_\chi$. The weights that appear are easily calculated and are shown in Figure \ref{fig:K}. The first term of the resolution is $\cT_{(2,-2)}$, so by deleting that term, we obtain a free resolution of the sheaf $\cK$. 

Note that the final term of the resolution is $\cT_{(0,-3)}$ (with multiplicity one), and this bundle does not appear in the other terms.

\begin{lem}\label{EndK} 
 $$\Hom(\cE_+, \cE_+)= \Hom(\cE_-, \cE_-)$$
and $\Ext^k(\cE_+, \cE_+) = 0 $ for $k>0$.
\end{lem}
\begin{proof} 
Write $\cS = q_*(L^{-4}) \otimes (\det T)^2$, so by definition we have a short exact sequence:
$$\cE_+\to \cE_- \to \cS $$
Our claim is that:
$$\mbox{(i)} \; \;\Ext^\bullet(\cE_+, \cS)=0 \aand \mbox{(ii)} \;\; \Ext^\bullet(\cS, \cE_-)=0$$
For (i), we claim more precisely that $\cS$ is right orthogonal to $\cT_\chi$ for every $\chi$ in $\nabla_+$ and every $\chi$ appearing in our free resolution of $\cK$. This holds by the first statement of \cite{SVdB17}*{Lem.~11.2.1} since every such $\chi$ satisfies the inequality:
$$\langle  \lambda_2, \chi\rangle\;> \; \langle \lambda_2, (2,-2)\rangle $$
Part (ii) is the claim that $\cS$ is left orthogonal to $\cT_\chi$ for every $\chi\in \nabla_-$. We use the fact \cite{W03} that $\cS^\vee$ is a shift of:
$$ q_*((T/L)^{3}) $$
Then left orthogonality to $\nabla_+$ follows as above, since $\nabla_+$ is invariant under 
$$\cT_{(a,b)}\mapsto \cT_{(a,b)}^\vee = \cT_{(-b,-a)}$$
and every $\chi\in \nabla_+$ satisfies $\langle  \lambda_2, \chi\rangle > \langle \lambda_2, (3,0)\rangle $. Finally, to get that $\Ext^\bullet(\cS, \cT_{(2,-2)})=0$ we need to use the rest of \cite{SVdB17}*{Lem.~11.2.1},  and observe that every weight in $(a,b) \in (V^{\lambda_2\geq 0})^\vee$ satisfies $a>0$. 
\end{proof}

\begin{lem}\label{GRRforK} $\cK$ lies in the grade restriction category $\cG_+$.
\end{lem}
\begin{proof}
It is clear from Figure \ref{fig:K} that $\cK$ satisfies the required grade restriction rules for $\lambda_0$ and $\lambda_2$ on both sides because all the terms in the free resolution lie in the correct region. But this is not true for $\lambda_1$; in fact, $\cK$ cannot satisfy the $\lambda_1$ grade restriction rule on the negative side. If it did, then the torsion sheaf $\cS$ would lie in $\cG_-$, which would be a contradiction since $\cS$ restricts to zero in $X_-$. 

However, for $\cK$ to lie in $\cG_+$, we only require that $\cS$ (and hence also $\cK$) satisfies the $\lambda_1$ grade restriction rule along the locally closed subset:
$$ Z_1^+ = V^{\lambda_1}\setminus \{\psi \in \Sigma_1\} $$
We verify this with an explicit computation. 

We are interested in the torsion sheaves between $\cO_{Z_1}$ and $\cS = q_*(L^{-4}) \otimes (\det T)^2$, and more specifically in the $\lambda_1$ weights that occur in them. Recall, in the notation of Section \ref{sec.GR}, that
$$ Z_1 = \set{ \notate{\star}{0}{\star}{0}{0}{0}} \aand Z_1^+ = \set{\notate{\star}{0}{\square}{0}{0}{0}}$$
Using the evident projective coordinates on the resolution \eqref{resolution}, we can write the map $q$ as:\footnote{The components of the 4-vector here are the coefficients of the cubic $(vt_1 -ut_2)(\alpha t_1^2 + \beta t_1t_2 + \gamma t_2^2)$ that vanishes on the line $u\!:\!v$.}
$$q: (x_1,..., x_6\,\vert\, \alpha, \beta, \gamma \,\vert\, u\!:\!v) \mapsto \mat{ ux_1 &... & ux_6 \\ vx_1 & ... & vx_6 }, \;
 (v\alpha, v\beta- u\alpha, v\gamma - u \beta, -u\gamma) $$
We care about the $\lambda_1$ grading, which on the source coordinates we can take to be
$$1, ..., 1 \,\vert\,0, -1, -2 \,\vert\, -\!1, 0 $$
(we have some freedom here because these are projective coordinates). Since we only care about $Z_1^+$, we can set $v=1$ and $\alpha\neq 0$, so $q$ becomes a map from $\C^9\times\C^*$ to $\C^{15}\times \C^*$. The line bundle $L$ is trivial on this locus, so $\cS$ is the pushdown of the free graded module $\cO(2)$. 
\pgap

To compute the torsion between $q_*\cO$ and $\cO_{Z_1}$, we take the Koszul resolution of the latter and pull it back via $q$. This produces the Koszul complex for the nine elements
$$x_1,..., x_6,  \beta- u\alpha,\, \gamma - u\beta, \, -u\gamma $$
in the graded ring $\C[x_1,..., x_6, \alpha^{\pm 1}, \beta, \gamma, v ] $. This is a regular sequence;  the vanishing locus is a (singular) curve so it has the expected dimension. Hence
$$\Tor^\bullet(q_*\cO, \cO_{Z_1^+}) = \Tor^0(q_*\cO, \cO_{Z_1^+})
=
\C[\gamma^{\pm 1}, u] / u^3 \; =\;  \cO_{Z_1}[u] / (x_1,..., x_6,  u^3 )$$
The $\lambda_1$ weight of the variable $u$ is $-1$ so the weights that occur here are exactly $0,-1,-2$. It follows that $\cS$ satisfies the $\lambda_1$ grade restriction rule. 
\end{proof}

\begin{rem}
As a check on this calculation, we can look instead at the locus $Z_1^-=Z_1\setminus\{\Psi=0\}$. This leads us to the Koszul complex for
$$vx_1,..., vx_6, v\beta- \alpha,\, v\gamma -  \beta,\, -\gamma $$
over the ring $\C[x_1,..., x_6, \alpha, \beta, \gamma, v ] $. This is not a regular sequence, and it has the same homology as the Koszul complex of $(vx_1,..., vx_6)$. The $\lambda_1$ weights here are $0,...,-5$ violating the grade restriction rule as expected.
\end{rem}

\subsection{Generation on \texorpdfstring{$X_+$}{X+}}

We now prove our final required fact, that the sheaf $E_+$ spans the category $D^b(X_+)$. 

\begin{prop}\label{E+spans} If $\cF\in D^b(X_+)$ satisfies
$$\Ext^\bullet(E_+, \cF) = 0$$
then $\cF\simeq 0$. 
\end{prop}
\begin{proof}
Recall that $X_+$ contains an open set
$$X_+^o \cong [\Hom(\C^6, \C^2) \,/\, \Gamma] $$
and the complement is the divisor $D$ \eqref{D}, which is a bundle over $\P^5$. 
We denote the two inclusions by:
\[
\begin{tikzcd}
	{X_+^o}  & {X_+} & D
		\arrow["f", hook, from=1-1, to=1-2]
	\arrow["g"', hook', from=1-3, to=1-2]
\end{tikzcd}
\]
The local cohomology of $\cF$ along $D$ fits into a triangle:
$$R\Gamma_D(\cF) \to \cF \to f_*f^*\cF $$
Now consider the subcategory $\langle E_+\rangle \subset D^b(X_+)$ generated by the summands of $E_+$. By assumption, $\cF$ is orthogonal to this subcategory. 
By Lemma \ref{O_D} below, $\langle E_+\rangle$ contains all objects that are pushed-forward from $D$, and hence it contains the whole subcategory $D^b_D(X_+)$ of objects supported on $D$. 
It follows that $R\Gamma_D(\cF) =0$.

The irreducible representations of $G$ indexed by $\nabla_+$, when restricted to $\Gamma$, provide us with all nine irreducible representations of $\Gamma$ (see Section \ref{repthy} and fact (iii) in particular). So $f^*E_+$ spans $D^b(X_+^o)$ and hence we also have $f_*f^*\cF=0$.  
\end{proof}

\begin{lem}\label{O_D}The image of the functor $g_*:D^b(D)\to  D^b(X_+)$ lies in the subcategory $\langle E_+\rangle$ generated by the summands of $E_+$.
\end{lem}

\begin{proof}
The divisor $D$ is cut out by a section of the line bundle $(\det T)^6 = \cT_{(6,6)}$ (Remark \ref{discriminant}), so we have a short exact sequence:
$$\cT_{(-6,-6)} \to \cO_{X_+} \to g_* \cO_D$$
Moreover, $D^b(D)$ is generated by the six line bundles $\cO_D,..., \cO_D(5)$. So we are reduced to showing that the line bundles 
$$\cT_{(k,k)} = \det(T)^k, \quad\mbox{for } k\in[-6,5] $$
lie in $\langle E_+\rangle$.

Recall the locus $GV^{\lambda_1 \leq 0}$ destabilised by $\lambda_1$ on the positive side. This is the set where $\psi$ has a triple root, \ie $\Hom(\C^6, T)\times\Sigma_1$. Ignoring the linear factor, its resolution is the line bundle $(T/L)^{-3}_{\P T}$. For every positive line bundle on the resolution, we get an exact sequence of vector bundles on $X_+$, of which the first three are:
$$ \cT_{(5,4)} \To \cT_{(4,2)} \To \cO $$
$$ \cT_{(5,5)} \To \cT_{(3,1)} \To \cT_{(1,0)}$$
\beq{seqs}\cT_{(5,3)} \to \cT_{(4,1)} \to \cT_{(2,0)}\eeq
(the first two are dual up to a twist). Starting from the set $\nabla_+$, we can use twists of these three sequences to recursively generate the bundles $\cT_\chi$ for the weights shown in Figure \ref{fig:gen} (L), which include $\cT_{(k,k)}$  for $k\in[-5,5]$.

\begin{figure}[ht]
\begin{tabular}{cc}

\begin{tikzpicture}[scale=0.75]
\fill[fill=gray,fill opacity=0.2]   (-3,-3) -- (3,3) -- (3,-3) -- cycle;

\draw[<->, thin, color=gray] (-3,0) -- (3,0);
\draw[<->, thin, color=gray] (0,-3) -- (0,3);
\draw[step=0.5, thin, color=gray!30!white] (-3,-3) grid (3,3);

\foreach \i in {  (-1,-1), (-0.5,-0.5), (0,0), (0.5,.5),(1,1), (-0.5,-1),(0,-0.5),(0.5,0),(1,0.5), (0,-1), (0.5,-0.5), (1,0), (0.5,-1), (1,-0.5), (1.5,1.5), (1.5,1), (1.5,0.5), (2,2), (2,1.5), (2.5,2.5), (-1.5,-1.5), (-1,-1.5), (-.5,-1.5), (-2,-2), (-1.5,-2), (-2.5,-2.5)}
{ \draw[very thick, red!90!black] \i circle (4 pt); }

\foreach \i in {  (-1,-1), (-0.5,-0.5), (0,0), (0.5,.5),(1,1), (-0.5,-1),(0,-0.5),(0.5,0),(1,0.5), (0,-1), (0.5,-0.5), (1,0), (0.5,-1), (1,-0.5)}
{ \filldraw[red!90!black] \i circle (3 pt); }

\end{tikzpicture}
& \hspace{1cm}
\begin{tikzpicture}[scale=0.75]
\fill[fill=gray,fill opacity=0.2]   (-3,-3) -- (3,3) -- (3,-3) -- cycle;

\draw[<->, thin, color=gray] (-3,0) -- (3,0);
\draw[<->, thin, color=gray] (0,-3) -- (0,3);
\draw[step=0.5, thin, color=gray!30!white] (-3,-3) grid (3,3);

\foreach \i in {  (-1,-1), (-0.5,-0.5), (0,0), (0.5,.5),(1,1), (-0.5,-1),(0,-0.5),(0.5,0),(1,0.5), (0,-1), (0.5,-0.5), (1,0), (0.5,-1), (1,-0.5), (1.5,1.5), (1.5,1), (1.5,0.5), (2,2), (2,1.5), (2.5,2.5), (-1.5,-1.5), (-1,-1.5), (-.5,-1.5), (-2,-2), (-1.5,-2), (-2.5,-2.5)}
{ \draw[very thick, red!90!black] \i circle (4 pt); }

\foreach \i in { (0,-1.5), (0,-1), (0,-.5), (0,0),(.5,.5), (1, .5),
(-1.5, -1.5), (-1,-1), (-.5, -1), (0,-1), (.5, -1)}
{ \filldraw[green!60!black] \i circle (2 pt); }

\filldraw[green!60!black] (0,-1.5) circle (4 pt);

\end{tikzpicture}
\end{tabular}

\caption{\footnotesize (L) Bundles that can be generated from $\nabla_+$ using the sequences \eqref{seqs} alone. (R) Generating $\cT_{(0,-3)}$ from the previous set and $\cK$.} \label{fig:gen}
\end{figure}
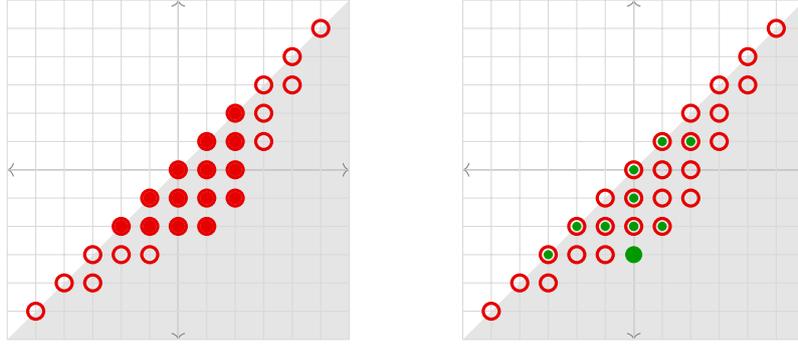

Now compare this with our free resolution of $\cK$ (Figure \ref{fig:K}). Since the final term is $\cT_{(0,-3)}$ we have a cone $[\cK\to \cT_{(0,-3)}]$ and this consists of bundles that we have already generated. So, by including $\cK$ as well, we can generate $\cT_{(0,-3)}$. See Figure \ref{fig:gen} (R). Then two more applications of the sequences \eqref{seqs} give us $\cT_{(-2,-4)}$ and finally $\cT_{(-6,-6)}$.
\end{proof}

\begin{rem} \label{selfExt} Consider the locus $GV^{\lambda_2\leq 0}$ destabilized by $\lambda_2$ on the positive side. By pushing down the right line bundle from the resolution of this locus, we can produce an exact sequence of bundles on $X_+$ which both starts and ends with $\cT_{(2,-2)}$. It represents a non-zero class in
$$\Ext^7_{X_+}\!\!\big(\cT_{(2,-2)}, \cT_{(2,-2)}\big) .$$
\end{rem}

\bibliography{bib_fano_of_lines}{}
\bibliographystyle{plain}

\end{document}